\documentclass{article}
\usepackage{color,amsmath,amssymb,hyperref,graphicx}
\numberwithin{equation}{section}

\author{Katharina Habermann, Lutz Habermann} 
\title{An Evolutionary Game-Theoretic Approach to Open Access}
\newcommand{\bx}{\textbf{\textit{x}}}
\newcommand{\by}{\textbf{\textit{y}}}
\newcommand{\R}{\mathbb{R}}
\newcommand{\rd}{\mathrm{d}}

\newcommand{\vp}{\varphi}
\newcommand{\om}{\omega}

\newtheorem{lm}{Lemma}[section]
\newtheorem{sz}[lm]{Proposition}

\newtheorem{fo}[lm]{Corollary}
\newtheorem{df}[lm]{Definition}
\newtheorem{bm0}[lm]{Remark}
\newtheorem{bs0}[lm]{Example}

\newenvironment{beweis}{{\em Proof.}}{\hspace*{\fill} $\square$}
\newenvironment{bm}{\begin{bm0}\rm}{\hspace*{\fill} $\square$ \end{bm0}}

\begin{document}
\maketitle
\begin{abstract}
The paper presents an evolutionary game-theoretic approach to open access
publishing as an asymmetric game between scientists and publishers. 
We show how the ordinary differential equations of the model presented can be
written as a system of Hamiltonian partial differential equations. The
understanding of the setting as a Hamiltonian system implies some properties 
reflecting the qualitative behavior of the system.
\end{abstract}
\section{Introduction}

The topic of open access publishing has been extensively and controversially
discussed. For general information see, for instance,
\cite{plattform,sparc,sparc-news}. There are many models of open access
named after different colors such as golden or green roads to open access
\cite{inphys,ceic}, but we do not want to go into more details, here.
Concerning scientific publishing particularly in mathematics, we refer to
\cite{birman,jackson} and \cite{navin} as well as to the references therein.

In \cite{oa-arxiv} the authors illustrate a game-theoretic approach to open
access publishing in order to understand different publication patterns within
different scientific disciplines. The underlying observation is that there are
communities where open access publishing is widely adopted, whereas other
scientific communities are far away from practicing any open access
publishing. First, different classical game settings are discussed
namely a zero sum game, a game similar to the Prisoners' Dilemma (up to sign), 
and a stag hunt game version, all describing a Nash equilibrium dilemma of
the non-open access communities. Second, these classical settings are
transferred into their quantum game extensions which allow to tackle the
dilemma mentioned that cannot be solved within the classical setting.

The classical approach in \cite{oa-arxiv} is formulated as a two-player game
where the players are authors (scientists) in symmetrical situations. The two
authors have the same set $\{s_1,s_2\}$ of strategies at their disposal. 
Consequently, the game looks exactly the same to both of its players.
For the open access game each player has to decide between the option to
publish open access as the first strategy $s_1$ or, as the second strategy
$s_2$, to conventionally publish in traditional journals where articles go
through peer reviewing. Depending on the strategy chosen by the co-player, the
players aim to maximize their success of the game. Due to the symmetric
situation both players do not only have the same set of strategies but also
have the same payoff matrices.

In contrast to the approach of a symmetric two-scientists game, the present
paper attempts to explain mathematically the open access play as a conflict of
interest between scientists and publishers, as it is emphasized by the
rapidly accumulating literature on open access and electronic publishing. 
Assuming that the open access problem is not a conflict between those scientists
who publish open access and those who do not, one takes the view that it is
more adequate to understand the problem of open access as an asymmetric
conflict between scientists and publishers.

The crucial point of this view is that it actually suggests a bimatrix game
describing this asymmetric setting.  

We do not explain why both players -- scientists and publishers -- behave as
they do in the game situation described, but we will give a mathematical
description of their interaction in the game. Moreover, we do not argue for
preferences in choosing one game strategy over another.
\section{The Game Setting}

We consider scientists and publishers as players in different positions them
having different strategy sets as well as different payoff matrices. 
Moreover, we take into consideration whole populations of players parts of
which choose either one or the other strategy.

The set of pure strategies for the scientists is $\{s_1,s_2\}$, where $s_1$ is 
publishing open access and $s_2$ stands for conventional
publishing. Furthermore, let $0\le x_1\le 1$ and $0\le x_2\le 1$ denote the
relative parts of the scientists' population playing strategy $s_1$ and $s_2$,
respectively. If, for example, half of all scientists publishes open
access and the other does not, one has
$\displaystyle{x_1=x_2=\frac{1}{2}}$. Or, if the relation is $20\%$ to $80\%$,
then $\displaystyle{x_1=\frac{1}{5}}$ and
$\displaystyle{x_2=\frac{4}{5}}$. Since $x_1$ and $x_2$ describe relative 
frequencies, we see that $x_1+x_2=1$ always holds true. Hence, the set of all
possible mixed strategies for the scientists' population is
$$ S=\left\{\bx={x_1 \choose x_2}\in\R^2\;:\; 0\le x_1,x_2\le 1 \mbox{ and }
  x_1+x_2=1\right\}. $$
The opponent population of publishers consists of two types playing pure
strategies $\{p_1,p_2\}$. Here strategy $p_1$ means the publisher accepts
or realizes open access publishing, whereas $p_2$ represents the strategy
of definite declining any way of open access. For the publishers' population let
$y_i$ be the frequency of strategy $p_i$, where $i=1,2$. Thus, all possible
mixed strategies for the publishers' population are given by the set
$$ P=\left\{\by={y_1 \choose y_2}\in\R^2\;:\; 0\le y_1,y_2\le 1 \mbox{ and }
  y_1+y_2=1\right\}\;. $$

In order to discuss the conflict of open access publishing, we proceed with
establishing the payoff matrices $A$ for the scientists and $B$ for the
publishers. Each payoff depends on the strategies chosen by the participating
players. A player in the scientists' population using strategy $s_i$ against a
player from the publishers' population using strategy $p_j$ obtains the payoff
$a_{ij}$, whereas the opponent obtains $b_{ji}$, with $i,j=1,2$.

As reasoned in \cite{oa-arxiv}, it is convenient to assume that the scientist
tries to maximize his scientific reputation. Let $R>0$ denote the reputation
payoff that the scientist is awarded. This payoff will be slightly reduced
by $0<r<R$ in case of publishing open access. Hence, we still have $R-r>0$ for
the total reputation.
In addition, we take into consideration the impact of an article by the
scientist. The impact gives some payoff $I>0$ both for the scientist and for
the publisher. Further, this impact obviously will reduce by some $0<\iota<I$,
if the concerning journal is not sufficiently available.
Moreover, we assume that $\iota<r$ saying that the reduced impact in case of
non-open access publishing is smaller than the loss of reputation in case of
open access publishing.
Open access publishing will cause expenses $L>0$ that are shared equally by both
players if both play their open access strategy.
Furthermore, $G>0$ expresses some moderate journal price scientists have
to pay for library subscriptions as well as the compensation for
expenses of the publisher. Finally, let us introduce $P>0$ for exorbitant
profit representing the prices of the most expensive journals coupled with
continuing non-open access publishing and taking the tremendously increasing
fees for access for granted. For our considerations, we make the realistic
assumption that $P$ is so large that the inequality
\begin{equation}\label{profit}
G+P-L>r-\iota
\end{equation}
holds true, which is equivalent to
\begin{equation*}
G+P>L+r-\iota\;.
\end{equation*}
We already have $r-\iota>0$ and the open access expenses $L$ actually
cannot exceed the upper bound given by this inequality.

If an open access publishing scientist meets a publisher accepting open
access, the scientist will get his reduced reputation, his 
research impact, and will pay for maintaining some open access publishing
system as well as for the journal subscription. The publisher takes the price
for the journal subscription and some amount given by the impact. Apart from
this, the publisher contributes to some extent to open 
access. Hence, the payoff is $\displaystyle{(R-r)+I-\frac{L}{2}-G}$ for the
scientist and $\displaystyle{G+I-\frac{L}{2}}$ for the publisher. 

If the scientist's strategy is open access whereas the publisher plays the
non-open access strategy, the scientist's payoff will be
$\displaystyle{(R-r)+I-L}$. The scientists have to run the open access
completely by themselves and do not have any expenses for journal
subscriptions. For the publisher nothing happens and hence the payoff is $0$.

A non-open access publishing scientist gets his full reputation $R$ but
only reduced research impact $I-\iota$. Additionally, the journal has to be
paid, which reduces the payoff by $G$.
If the publisher's counter-strategy is open access, the payoff for the scientist
will be $\displaystyle{R+(I-\iota)-G}$. In this case the payoff for the
publisher will be $\displaystyle{G+(I-\iota)-L}$ paying the total sum for
open access acceptance.
Moreover, the payoff is $\displaystyle{R+(I-\iota)-G-P}$ for the scientist and
$\displaystyle{G+(I-\iota)+P}$ for the publisher in case of the non-open
access strategy of the publisher. 

Altogether, we have the complete income-and-loss statement as given in
Table~\ref{tabelle1}. 

\begin{table}[h]
\caption{\label{tabelle1} {\sf Complete payoffs}}
\begin{center}
\begin{tabular}{|c||c|c|}
\hline
& payoff & payoff \\
strategies & \quad for the scientist \quad & \quad for the publisher \quad  \\ 
\hline\hline && \\
$s_1\longleftrightarrow p_1$ & $\displaystyle{(R-r)+I-\frac{L}{2}-G}$ &
$\displaystyle{G+I-\frac{L}{2}}$ \\ && \\
\hline && \\
$s_1\longleftrightarrow p_2$ & $\displaystyle{(R-r)+I-L}$ &
$\displaystyle{0}$ \\ && \\
\hline && \\
$s_2\longleftrightarrow p_1$ & $\displaystyle{R+(I-\iota)-G}$ &
$\displaystyle{G+(I-\iota)-L}$ \\ && \\
\hline && \\
$s_2\longleftrightarrow p_2$ & $\displaystyle{R+(I-\iota)-G-P}$ &
$\displaystyle{G+(I-\iota)+P}$ \\ && \\
\hline
\end{tabular}
\end{center}
\end{table}

Hence, the corresponding payoff matrices are
$$ A=\left(\begin{array}{cc} a_{11} & a_{12} \\ a_{21} & a_{22} \end{array}
   \right)=\left(\begin{array}{ccc} \displaystyle{(R-r)+I-\frac{L}{2}-G} & &
   \displaystyle{(R-r)+I-L} \\ && \\ \displaystyle{R+(I-\iota)-G} & &
   \displaystyle{R+(I-\iota)-G-P}\end{array}\right) $$
and
$$ B=\left(\begin{array}{cc} b_{11} & b_{12} \\ b_{21} & b_{22} \end{array}
   \right)=\left(\begin{array}{ccc} \displaystyle{G+I-\frac{L}{2}} & & 
   \displaystyle{G+(I-\iota)-L} \\ && \\ \displaystyle{0} & &
   \displaystyle{G+(I-\iota)+P}\end{array}\right)\;. $$

Now, both payoff matrices determine the payoffs for the entire populations.
For a general description of how this happens we warmly recommend the book
\cite{hofb-sieg} by Hofbauer and Sigmund. 
The corresponding strategies are spread within both populations in accordance
with the frequencies of the strategies. Hence, if the scientists' population is
in state $\bx\in S$ and the publishers' population is in state $\by\in P$, then
the payoff for the entire population of scientists will be
$$ \langle\bx,A\by\rangle=\bx^TA\by $$
and that for the publishers' population will be
$$ \langle\by,B\bx\rangle=\by^TB\bx=\bx^TB^T\by\;. $$

As already mentioned, players will choose their game strategies with the
intention of maximizing their average payoff. This would be easy if one knew
the strategy that the opponent player is going to choose. In case both
players of the game choose simultaneously strategies of best reply to the
choices of the others, a pair of strategies occurs where both players are
incited to keep going into the same direction. Such a pair of strategy 
choices is known as Nash equilibrium. This is technically expressed by the
following definition.
\begin{df}\label{def-nash}
A {\it Nash equilibrium} is a pair $(\overline{\bx},\overline{\by})\in S\times
P$ such that 
$$ \bx^TA\overline{\by}\le\overline{\bx}^TA\overline{\by} \quad\mbox{ for all }
\bx\in S $$
as well as
$$ \by^TB\overline{\bx}\le\overline{\by}^TB\overline{\bx} \quad\mbox{ for
  all } \by\in P $$
is fulfilled.
\end{df}

The following considerations aim to give a Nash equilibrium for our open
access game setting.

\begin{lm}\label{lemma-positive}
Each of the differences $a_{12}-a_{22}$, $a_{21}-a_{11}$, $b_{22}-b_{21}$, and
$b_{11}-b_{12}$ is positive.
\end{lm}
\begin{beweis}
Indeed,
$$ a_{12}-a_{22}=-r-L+\iota+G+P>0 $$
by (\ref{profit}) and
$$ a_{21}-a_{11}=-\iota+r+\frac{L}{2}>0 $$
by $r>\iota$.
Furthermore,
$$ b_{22}-b_{21}=G+I-\iota+P $$
which is positive since $I>\iota$.
Finally,
$$ b_{11}-b_{12}=\frac{L}{2}+\iota $$
which is obviously positive.
\end{beweis}

\begin{lm}\label{nash-equ}
Let $(\overline{\bx},\overline{\by})\in S\times P$ be given by
$$ \overline{\bx}={x_0 \choose 1-x_0} \quad\mbox{ and }\quad 
   \overline{\by}={y_0 \choose 1-y_0}\;, $$
where
$$ x_0=\frac{b_{22}-b_{12}}{b_{22}-b_{12}+b_{11}-b_{21}} \quad\mbox{ and }\quad
   y_0=\frac{a_{12}-a_{22}}{a_{12}-a_{22}+a_{21}-a_{11}}\;. $$
Then
$$ \bx^TA\overline{\by}=\overline{\bx}^TA\overline{\by} \quad\mbox{ for all }
\bx\in S $$
and
$$ \by^TB\overline{\bx}=\overline{\by}^TB\overline{\bx} \quad\mbox{ for
  all } \by\in P $$
hold true.
\end{lm}
\begin{beweis}
First, by Lemma~\ref{lemma-positive} both denominators of the fractions
defining $x_0$ and $y_0$ are positive.
Then we have
\begin{eqnarray*}
A\overline{\by}
&=& \left(\begin{array}{cc} a_{11} & a_{12} \\ a_{21} & a_{22} \end{array}
     \right){y_0 \choose 1-y_0} \\
&=& {y_0(a_{11}-a_{12})+a_{12} \choose y_0(a_{21}-a_{22})+a_{22}} 
  \quad=\quad \frac{\det A}{a_{12}-a_{22}+a_{21}-a_{11}}{1 \choose 1} 
\end{eqnarray*}
as well as
\begin{eqnarray*}
B\overline{\bx}
&=& \left(\begin{array}{cc} b_{11} & b_{12} \\ b_{21} & b_{22} \end{array}
     \right){x_0 \choose 1-x_0} \\
&=& {x_0(b_{11}-b_{12})+b_{12} \choose x_0(b_{21}-b_{22})+b_{22}}
   \quad=\quad \frac{-\det B}{b_{22}-b_{12}+b_{11}-b_{21}}{1 \choose 1}\;.
\end{eqnarray*}
Hence,
$$ \bx^TA\overline{\by} 
   = \frac{\det A}{a_{12}-a_{22}+a_{21}-a_{11}}(x+(1-x))
   = \frac{\det A}{a_{12}-a_{22}+a_{21}-a_{11}} $$
does not depend on $\bx=\displaystyle{x \choose 1-x}$ and
$$ \by^TB\overline{\bx}
   = \frac{\det B}{b_{12}-b_{22}+b_{21}-b_{11}}(y+(1-y))
   = \frac{\det B}{b_{12}-b_{22}+b_{21}-b_{11}} $$
is independent of $\by=\displaystyle{y \choose 1-y}$.
\end{beweis}

\begin{fo}
The pair $(\overline{\bx},\overline{\by})\in S\times P$ given in
Lemma~\ref{nash-equ} is a Nash equilibrium.
\hfill $\square$
\end{fo}

\begin{bm}
For this game there is no strict, i.e. given by pure strategies, Nash
equilibrium. The Nash equilibrium $(\overline{\bx},\overline{\by})\in S\times
P$ specified in Lemma~\ref{nash-equ} is a unique mixed Nash equilibrium and
determined by the following equations, c.f. \cite{hofb-sieg}, Chapter~10.2
(Page~116, above),
\begin{eqnarray*}
a_{11}y_0+a_{12}(1-y_0) &=& a_{21}y_0+a_{22}(1-y_0) \\
b_{11}x_0+b_{12}(1-x_0) &=& b_{21}x_0+b_{22}(1-x_0)\;.
\end{eqnarray*}
\end{bm}
\section{The Dynamics of the Game} 

Due to the payoffs of the game, there is an intrinsic dynamics within the
system. Actually, starting for instance with the situation that the
scientist has chosen strategy $s_1$ whereas the publisher's strategy is
$p_1$, the outcome is 
$$ (R-r)+I-\frac{L}{2}-G $$
for the scientist and  
$$ G+I-\frac{L}{2} $$
for the publisher.
In this situation the scientist could increase his payoff by changing
his strategy to $s_2$, since $r-\iota>0$ and $L>0$.
Indeed, the inequality
$$ R-r+I-\frac{L}{2}-G < R+I-\iota-G $$
is equivalent to 
$$ -\frac{L}{2} < r-\iota\;, $$
which holds true obviously for $r-\iota>0$ and $L>0$.
After this change in the scientist's strategy the payoff is
$$ R+I-\iota-G $$
for the scientist. Holding strategy $p_1$ the publisher's payoff is
$$ G+(I-\iota)-L\;. $$ 
In this situation, the publisher sees the chance to change strategy to
$p_2$. The reason is that the publisher's payoff of $G+(I-\iota)-L$ in the game
$s_2\leftrightarrow p_1$ can be increased to $G+(I-\iota)+P$ in a game
$s_2\leftrightarrow p_2$. In fact, $L>0$ and $P>0$ give
$$ G+(I-\iota)-L<G+(I-\iota)+P\;. $$
However, finding oneself in the situation of the game $s_2\leftrightarrow p_2$
the scientist notices that $G$ and $P$ are unnecessary costs, $L$ would not be
so huge and the loss of reputation will be moderate. Hence, the scientist will
prefer to choose strategy $s_1$. Altogether, changing strategies in order to
maximize outcome ends up moving in a circle. 

The dynamics for pairs of pure strategies looks like

\begin{center}
\begin{tabular}{ccc}
$(s_1,p_1)$ & $\leftarrow$ & $(s_1,p_2)$ \\
$\downarrow$ & & $\uparrow$ \\
$(s_2,p_1)$ & $\rightarrow$ & $(s_2,p_2)$ 
\end{tabular}
\end{center}

where the horizontal arrows picture the change of the publishers' strategies and
the vertical arrows illustrate the deviation in the scientists' strategies.

This dealing with an ``oscillating'' system requires a non-static approach.
Therefore, we consider the frequencies of the strategies within both populations
as time-dependent quantities which are differentiable functions of $t\in\R$,
i.e. $x_i=x_i(t)$ and $y_j=y_j(t)$ for $i,j=1,2$.
The first derivatives of these functions
$$ \dot x_i=\frac{\rd x_i}{\rd t} \quad\mbox{ and }\quad
   \dot y_j=\frac{\rd y_j}{\rd t} \quad\mbox{ for }i,j=1,2 $$
describe the rate of growth of the respective frequencies. 

Modeling the dynamics and the cyclic behavior of that system, we follow the
ideas of Hofbauer and Sigmund developed in \cite{hofb-sieg}, Chapter~10.3.
The increase of strategy $s_i$ of the scientists' population is given by
the per capita growth rate $\displaystyle{\frac{\dot x_i}{x_i}}$ and equals
the difference between its payoff $(A\by)_i$ and the average payoff
$\bx^TA\by$. The same holds true for the publishers' population.
This leads to the system of ordinary differential equations
\begin{eqnarray*}
\dot x_i &=& x_i\left((A\by)_i-\bx^TA\by\right) \\
\dot y_j &=& y_j\left((B^T\bx)_j-\by^TB^T\bx\right)
\end{eqnarray*}
for $i,j=1,2$. 

Let us look closer at these differential equations. We have
$$ A\by=\left(\begin{array}{cc} a_{11} & a_{12} \\ a_{21} & a_{22} \end{array}
  \right){y_1 \choose y_2}={a_{11}y_1+a_{12}y_2 \choose a_{21}y_1+a_{22}y_2} $$
and
$$ B^T\bx=\left(\begin{array}{cc} b_{11} & b_{21} \\ b_{12} & b_{22} 
  \end{array}\right){x_1 \choose x_2}={b_{11}x_1+b_{21}x_2 \choose b_{12}x_1
  +b_{22}x_2}\;. $$
Hence,
$$ (A\by)_i-\bx^TA\by
   =y_1(a_{i1}-x_1a_{11}-x_2a_{21})+y_2(a_{i2}-x_1a_{12}-x_2a_{22}) $$
and
$$ (B^T\bx)_j-\by^TB^T\bx
   =x_1(b_{1j}-y_1b_{11}-y_2b_{12})+x_2(b_{2j}-y_1b_{21}-y_2b_{22})\;. $$
Since $x_1+x_2=1$ and $y_1+y_2=1$, it is admissible to introduce new variables
$x:=x_1$ and $y:=y_1$ (then $x_2=1-x$ and $y_2=1-y$) and consider the
corresponding differential equations for $x$ and $y$.
These are
\begin{eqnarray*}
\dot x
&=& \dot x_1 \\
&=& x_1\left((A\by)_1-\bx^TA\by\right) \\
&=& x\left(y_1(a_{11}-x_1a_{11}-x_2a_{21})
    +y_2(a_{12}-x_1a_{12}-x_2a_{22})\right) \\
&=& x\left(y(a_{11}-xa_{11}-(1-x)a_{21})
    +(1-y)(a_{12}-xa_{12}-(1-x)a_{22})\right) \\
&=& x(1-x)\left(y(a_{11}-a_{21})+(1-y)(a_{12}-a_{22})\right) \\
&=& x(1-x)\left((a_{12}-a_{22})-(a_{12}-a_{22}+a_{21}-a_{11})y\right)
\end{eqnarray*}
and similarly
$$ \dot y \quad=\quad \dot y_1 \quad=\quad 
   y(1-y)\left(-(b_{22}-b_{21})+(b_{22}-b_{21}+b_{11}-b_{12})x\right)\;. $$
Now, let us introduce new constants by
\begin{eqnarray*}
a &:=& a_{12}-a_{22} \\
b &:=& a_{21}-a_{11} \\
c &:=& b_{22}-b_{21} \\
d &:=& b_{11}-b_{12}\;,
\end{eqnarray*}
which are positive by Lemma~\ref{lemma-positive}.
Hence, only two equations are remaining namely
\begin{eqnarray}
\dot x &=& x(1-x)(a-(a+b)y) \label{sys-1} \\
\dot y &=& y(1-y)(-c+(c+d)x)\;. \label{sys-2}
\end{eqnarray}
This is an adequate form for investigating the qualitative behavior of
the system. 

\medskip

The subsequent section is devoted to a further investigation of the dynamics
of the system given by Equations~(\ref{sys-1})~and~(\ref{sys-2}).

\section{Symplectic Reformulation}

We consider the square $Q=\{(x,y)\in\R^2\;:\;0\le x,y\le 1\}$. Obviously,
$Q\cong S\times P$ by 
$$ (x,y)\in Q \quad\mapsto\quad \left({x \choose 1-x},{y \choose 1-y}\right)
   \in S\times P\;. $$

We use symplectic techniques for investigating what happens in the
interior ${\rm int}Q$ of $Q$.
More precisely, in order to prove that the orbits of the
System~(\ref{sys-1})~and~(\ref{sys-2}) are closed and do not leave 
${\rm int}Q$ we describe the whole setting as a Hamiltonian system. 

The interior of $Q$ is the open subset
$$ M=\{(x,y)\in\R^2\;:\;0<x,y<1\} $$
of $\R^2$.
The function $\vp:M\to\R$ given by $\vp(x,y)=xy(1-x)(1-y)$ is positive. Hence,
$$ \om=\frac{1}{\vp}\rd x\wedge\rd y $$
defines a $2$-form $\om$ on $M$.

Since we are in dimension $2$ the following is clear.
\begin{lm}
The pair $(M,\om)$ is a symplectic manifold. \hfill $\square$
\end{lm}

We are going to derive that Equations~(\ref{sys-1})~and~(\ref{sys-2}) are the
equations of motion of the Hamiltonian 
system $(M,\om,H)$ where the Hamiltonian $H:M\to\R$ is defined as
$$ H(x,y)=c\ln(x)+d\ln(1-x)+a\ln(y)+b\ln(1-y) 
   \qquad\mbox{ for }(x,y)\in M\;. $$

The notion of Hamiltonian systems is well-known and widely used in
mathematical physics to describe the behavior of systems in classical
mechanics. In physics the Hamiltonian $H$ describes the energy of the system
under consideration. Hence, it would be interesting to understand,
what intrinsic property is determined by the ``energy'' of open access
publishing.

\begin{lm}\label{ham_vf}
The Hamiltonian vector field of $(M,\om,H)$ is given by
$$ X_H=\left(\begin{array}{r}\displaystyle{\vp\frac{\partial H}{\partial y}} 
       \\ \\ \displaystyle{-\vp\frac{\partial H}{\partial x}} \end{array}
       \right)\;. $$
\end{lm}
\begin{beweis}
Generally, $X_H$ is defined by the equation
$$ \om(X_H,\quad)=\rd H\;. $$
Clearly,
$$ \rd H=\frac{\partial H}{\partial x}\rd x+\frac{\partial H}{\partial y}\rd
y\;. $$
If $V$ is any vector field on $M$ with functions $v_1$ and $v_2$
as components, then
$$ \om(V,\quad)=\frac{1}{\vp}\rd x\wedge\rd y(V,\quad)
   =\frac{1}{\vp}\{v_1\rd y-v_2\rd x\}\;. $$
Hence,
$$ \om(V,\quad)=\rd H $$
if and only if
$$ \frac{\partial H}{\partial x}=-\frac{1}{\vp}v_2 \quad\mbox{ and }\quad
   \frac{\partial H}{\partial y}=\frac{1}{\vp}v_1\;, $$
which proves the assertion.
\end{beweis}

Let us remind that motions of the Hamiltonian System $(M,\om,H)$ are curves
$\gamma(t)$ in $M$ satisfying the differential equation
$$ \frac{\rd}{\rd t}\gamma(t)_{|t=s}=X_H(\gamma(s))\;. $$

Altogether, we now can prove the following proposition.
\begin{sz}
The solutions of (\ref{sys-1})~and~(\ref{sys-2}) are exactly the motions of
the Hamiltonian system $(M,\om,H)$.
\end{sz}
\begin{beweis}
We have
\begin{eqnarray*}
\vp\frac{\partial H}{\partial y}
&=& xy(1-x)(1-y)\left\{\frac{a}{y}-\frac{b}{1-y}\right\} \\
&=& ax(1-x)(1-y)-bxy(1-x) \\
&=& x(1-x)(a-ay-by)
\end{eqnarray*}
and
\begin{eqnarray*}
-\vp\frac{\partial H}{\partial x}
&=& -xy(1-x)(1-y)\left\{\frac{c}{x}-\frac{d}{1-x}\right\} \\
&=& -cy(1-x)(1-y)+dxy(1-y) \\
&=& y(1-y)(-c+cx+dx)\;.
\end{eqnarray*}
Thus, Equations~(\ref{sys-1})~and~(\ref{sys-2}) are equivalent to
$$ \vp\frac{\partial H}{\partial y}=\dot x \quad\mbox{ and }\quad
   -\vp\frac{\partial H}{\partial x}=\dot y\;, $$
which implies the statement by Lemma~\ref{ham_vf}.
\end{beweis}

As a consequence, the solutions in $M$ of the
System~(\ref{sys-1})~and~(\ref{sys-2}) correspond to level set curves of the
Hamiltonian $H$, i.e. the Hamiltonian $H$ is constant along solutions.
\begin{fo}
Let $\left(x(t),y(t)\right)$ give a solution of
Equations~(\ref{sys-1})~and~(\ref{sys-2}). Then 
$$ \frac{\rd}{\rd t}H(x(t),y(t))=0 \;. $$ 
\end{fo}
\begin{beweis}
Indeed,
\begin{eqnarray*}
\lefteqn{\frac{\rd}{\rd t}H(x(t),y(t))=} \\
&=& \frac{\rd}{\rd t}\left\{c\ln(x(t))+d\ln(1-x(t))+a\ln(y(t))+b\ln(1-y(t))
    \right\} \\
&=& c\frac{\dot x}{x}-d\frac{\dot x}{1-x}+a\frac{\dot y}{y}
    -b\frac{\dot y}{1-y} \\
&=& (c(1-x)-dx)(a-(a+b)y)+(a(1-y)-by)(-c+(c+d)x) \\
&=& (c-cx-dx)(a-ay-by)+(a-ay-by)(-c+cx+dx) \\
&=& 0
\end{eqnarray*}
by using (\ref{sys-1}) as well as (\ref{sys-2}) for replacing $\dot x$ and
$\dot y$.
\end{beweis}

In the remainder of this section we want to discuss the shape of the level
set curves of $H$. In particular, we will show that these curves are closed,
which implies that the solutions of (\ref{sys-1})~and~(\ref{sys-2}) are
periodic.

First, we are going to determine the extremal points of $H$.
By 
$$ \frac{\partial H}{\partial y}=\frac{a}{y}-\frac{b}{1-y} $$
and
$$ \frac{\partial H}{\partial x}=\frac{c}{x}-\frac{d}{1-x} $$
we obtain
$$ \frac{a}{y_0}=\frac{b}{1-y_0} \quad\mbox{ and }\quad
   \frac{c}{x_0}=\frac{d}{1-x_0} $$
as necessary conditions for a critical point $(x_0,y_0)$.
This leads to
$$ y_0=\frac{a}{a+b} \quad\mbox{ and }\quad x_0=\frac{c}{c+d} $$
which corresponds exactly to our Nash equilibrium described in
Lemma~\ref{nash-equ}. By
\begin{eqnarray*}
\frac{\partial^2 H}{\partial^2 x}
&=& -\frac{c}{x^2}-\frac{d}{(1-x)^2}\;, \\
\frac{\partial^2 H}{\partial^2 y}
&=& -\frac{a}{y^2}-\frac{b}{(1-y)^2}\;, \quad\mbox{ and } \\
\frac{\partial^2 H}{\partial x\partial y}
\quad=\quad
\frac{\partial^2 H}{\partial y\partial x}
&=& 0\;,
\end{eqnarray*}
we have 
\begin{eqnarray*}
\frac{\partial^2 H}{\partial^2 x}(x_0,y_0) &=& -\frac{(c+d)^3}{cd}\quad<\quad 0
\quad\mbox{ and } \\
\frac{\partial^2 H}{\partial^2 y}(x_0,y_0) &=& -\frac{(a+b)^3}{ab}\quad<\quad 0
\end{eqnarray*}
at the point $(x_0,y_0)$.
This implies that the matrix
$$ {\rm
  Hess}H_{(x_0,y_0)}=\left(\begin{array}{cc}\displaystyle{-\frac{(c+d)^3}{cd}}
      & 0 \\ 0 & \displaystyle{-\frac{(a+b)^3}{ab}}\end{array}\right) $$
is negative definite showing that $H$ takes in $(x_0,y_0)$ a maximum as the
unique extremum. Second, since
$$ \lim_{t\to 0+}\ln(t)=-\infty\, $$
we have 
$$ \lim_{(x,y)\to q}H(x,y)=-\infty $$
for each $q$ in the boundary $\partial Q$ of $Q$. 
Hence, the graph of the Hamiltonian $H$ looks like a hill over the square
$Q$.

\begin{figure}[h]
\begin{center}
\includegraphics[width=70mm]{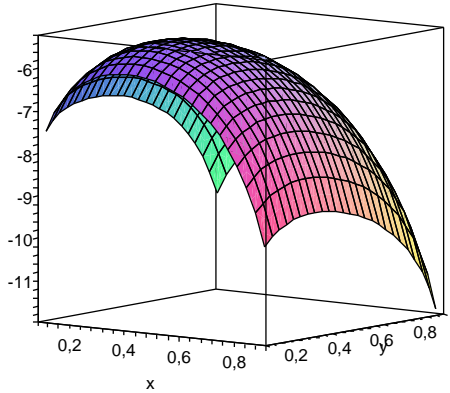}
\caption{\small\sf Graph of the Hamiltonian $H$ with parameters $a=1$, $b=2$,
  $c=2$, and $d=3$}\label{fig3}  
\end{center}
\end{figure}

The top of the hill represents the unique maximum. 
Towards the boundary $\partial Q$, the hill crashes into bottomless depths.

\begin{figure}[h]
\begin{center}
\includegraphics[width=60mm]{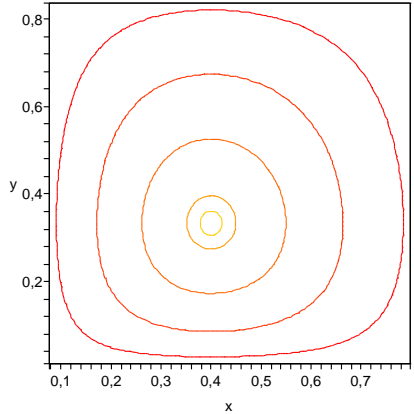}
\caption{\small\sf Periodic orbits surrounding the Nash equilibrium
  $(x_0,y_0)$ with parameters $a=1$, $b=2$, $c=2$, and $d=3$}\label{fig4} 
\end{center}
\end{figure}

Since all solutions of our system correspond to level set curves, they are
visualized as contour lines of the hill described.
Hence, all orbits surround the point $(x_0,y_0)$ and remain in $M$. 

\section{Further Discussion}

Obviously, the model presented here is a basic approach to discuss the
phenomenon of multiple arrangements of realizing open access publishing in
different scientific communities. 
However, in contrast to the approach presented in \cite{oa-arxiv} it treats open
access publishing in a canonical way as a game between scientists and
publishers. 
The behavior of oscillation on closed orbits is a rather rough and
simplified view of the whole scenario. We do not consider this model to be
universally valid, but throughout a certain period of time it seems to 
allow a reasonably satisfying description of insights which are intuitively
clear. 

This model can be modified.
For instance, taking into account that both the parameter $r$ reducing the
scientist's reputation in case of open access publishing and the value of
$\iota$ expressing the reduced impact in case of non-open access publishing
may change in the context of the game.
For example, if the rate of open access publishing scientists increases, $r$
will be supposed to decrease etc. Playing a little bit and slightly
altering the model presented here may result in further refinements.

Moreover, the huge number of evolutionary game-theoretic approaches in the
literature allows a lot of modifications and further developments of these
concepts. They vary from biological game theory to concepts that are more 
applicable in economic contexts. 
These models are good sources and suitable to establish further
improvements in a more fundamental manner, see e.g. \cite{oeco,hof-si} and
\cite{crash} as well as the references therein.
\subsection*{Acknowledgment}
The authors have benefited from many discussions with both scientists and
information professionals and particularly from the working context of the
first author as a mathematician as well as a subject librarian. 
The idea developed in this note arose during
the 5th spring school on Geometry and Mathematical Physics at Hiddensee in May 
2007. Both authors thank the Biological Station of Hiddensee for the kind
hospitality as well as for the inspiring working atmosphere, the organizer of
the workshop Stephan Block for his invitation, and the participants for the
intensive discussion.


\vfill

{\sc G\"ottingen State and University Library, 
Platz der G\"ottinger Sieben 1,  
37073 G\"ottingen,
Germany} \\
{\sf e-mail: habermann@sub.uni-goettingen.de}

{\sc Department of Mathematics,
University of Hannover,
Welfengarten 1,
30167 Hannover,
Germany} \\ 
{\sf e-mail: habermann@math.uni-hannover.de}

\end{document}